\documentclass[12pt, reqno]{amsart}

\usepackage[margin=1in]{geometry}

\usepackage[english]{babel}
\usepackage[utf8x]{inputenc}
\usepackage{amsmath}
\usepackage{amssymb}
\usepackage{amsthm}
\usepackage{enumerate}
\usepackage{hyperref}
\usepackage[numbers]{natbib}

\newtheorem{theorem}{Theorem}
\newtheorem*{theorem*}{Theorem}

\newtheorem*{corollary*}{Corollary}

\newtheorem{lemma}[theorem]{Lemma}

\newcommand{\CR}{\mathcal{R}}
\newcommand{\BH}{\mathbb{H}}

\newcommand{\BZ}{\mathbb{Z}}

\renewcommand{\Re}{\,\mathrm{Re}\,}

\newcommand{\symtwo}{\operatorname{sym}^{2}}

\allowdisplaybreaks[1]

\title[Restriction of Hecke eigenforms to horocycles]{Restriction of Hecke eigenforms to horocycles}

\author{Ho Chung Siu and Kannan Soundararajan}
\address{Department of Mathematics, Stanford University, Stanford, CA 94305}
\email{soarersiuhc@gmail.com}
\address{Department of Mathematics, Stanford University, Stanford, CA 94305}
\email{ksound@stanford.edu}
\thanks{The second author is supported in part by a grant from the National Science Foundation, and 
a Simons Investigator award from the Simons Foundation}
\date{\today}

\begin{document}
\maketitle
\begin{abstract}
    We prove a sharp upper bound on the $L^2$-norm of Hecke eigenforms restricted to a horocycle, as the weight tends to infinity.
\end{abstract}

\section{Introduction}

\noindent A central problem in ``quantum chaos" is to understand the limiting behavior of eigenfunctions. 
An important example that has attracted a lot of attention 
is that of Maass cusp forms with large Laplace eigenvalue on the modular surface 
$X = SL_2(\BZ) \backslash \BH$.  Let $\phi$ denote such a Maass form, with eigenvalue $\lambda$, and normalized 
to have $L^2$-norm $1$: that is, $\int_X |\phi(z)|^2 \frac{dxdy}{y^2} = 1$.   Then the Quantum Unique Ergodicity (QUE) conjecture of Rudnick and Sarnak \cite{RudSar} 
states that the measure $\mu_\phi = |\phi(z)|^2 \frac{dx dy}{y^2}$ tends to the uniform measure on $X$ as $\lambda \to \infty$.  If $\phi$ is also 
assumed to be an eigenfunction of all the Hecke operators, then QUE holds by the work of Lindenstrauss \cite{Lin}, with a 
final step on escape of mass provided by Soundararajan \cite{Sou}.   Thus, the measure $\mu_\phi$ does not concentrate on subsets of $X$ with 
small measure, but is uniformly spread out.  A finer problem is to understand how much the measure can concentrate on sub-manifolds; for example, 
on a geodesic, or a closed horocyle, or even at just a point (that is, bounding the $L^{\infty}$ norm).  The letter of Sarnak to Reznikov \cite{Sar} draws 
attention to such restriction problems, and these problems (and generalizations) have been studied extensively in recent years, see for example \cite{BoRu}, \cite{BurGerTsv}, 
\cite{GRS}, \cite{LiLiuYou}, \cite{LiYou}, \cite{TotZel_1}, \cite{TotZel_2}, \cite{You}.


This note is concerned with a related question for holomorphic modular forms for $SL_2({\Bbb Z})$ that are also eigenfunctions of 
all Hecke operators, when the weight $k$ becomes large.  Let $f$ be a Hecke eigenform of weight $k$ on the modular surface $X$, with $L^2$-norm $1$: that is, 
$$ 
\int_{X} y^k |f(z)|^2 \frac{dx dy}{y^2} = 1.
$$ 
To $f$, we associate the measure $\mu_f = y^k |f(z)|^2 \frac{dxdy}{y^2}$.  The analog here of QUE 
states that $\mu_f$ tends to the uniform measure $\frac{3}{\pi} \frac{dx dy}{y^2}$ as $k\to \infty$, and this is known 
to hold by the work of  Holowinsky and Soundararajan \cite{HolSou}.  As with Maass forms, one may now ask for finer restriction 
theorems for holomorphic Hecke eigenforms.  We study the problem of bounding the $L^2$-norm of Hecke eigenforms on 
a fixed horocycle, and establish the following uniform bound.

\begin{theorem}
\label{thm:main-thm}
Let $f$ be a Hecke eigenform of weight $k$ on $X = SL_2(\BZ) \backslash \BH$ with $L^2$-norm normalized to be $1$.  Let $\delta >0$ be 
fixed.  Uniformly in the range $1/k \le y\le k^{1/2-\delta}$ we have 
$$
\int_0^1 y^k| f(z)|^2 dx \le C(\delta),
$$ 
for some constant $C(\delta)$.  
\end{theorem}

Our result gives a uniform bound for the $L^2$-norm restricted to horocycles, answering a question from Sarnak \cite{Sar}.  
In the Maass form situation, Ghosh, Reznikov and Sarnak \cite{GRS} establish weaker restriction bounds (of size $\lambda^{\epsilon}$) 
for the corresponding problem, and Sarnak \cite{Sar} notes that uniform boundedness there follows from the Ramanujan conjecture and 
a sub-convexity bound (in eigenvalue aspect) for the Rankin-Selberg $L$-function $L(s,\phi\times \phi)$.  One might hope to strengthen and extend Theorem 1 in the following two ways.  First, Young \cite[Conjecture 1.4]{You} has conjectured that for any fixed $y > 0$, the restriction of $\mu_f$ to the horocycle $[0,1] + iy$ still tends to the uniform measure, as $k \to \infty$: in particular,  as $k\to \infty$
 $$
 \int_0^1 y^k | f(z)|^2 dx \to \frac{3}{\pi}. 
 $$
Second, one might expect that two different eigenforms $f$ and $g$ of weight $k$ are approximately orthogonal on the horocycle $[0,1]+iy$, so that (as $k \to \infty$) 
$$ 
\int_0^1 y^k f(x+iy) \overline{g(x+iy)} dx \to 0. 
$$ 
Our proof, which relies crucially on bounds for mean-values of non-negative multiplicative functions in short intervals, does not allow us to 
address these refined conjectures.

\section{Preliminaries}

\noindent Let $f$ be a Hecke eigenform of weight $k$ on $X = SL_2(\BZ) \backslash \BH$. Write
$$
L(s,f) = \sum_{n=1}^{\infty} \frac{\lambda_f(n)}{n^s} = \prod_p \Big(1 - \frac{\alpha_p}{p^s}\Big)^{-1}\Big(1 - \frac{\beta_p}{p^s}\Big)^{-1},
$$
where $\lambda_f(n)$ are the Hecke eigenvalues for $f$, and $\alpha_p$, $\beta_p = \alpha_p^{-1}$ are the Satake parameters. 
Our $L$-function has been normalized such that the Deligne bound reads $|\lambda_f(n)| \le d(n)$ (the divisor function), or equivalently  that $|\alpha_p| = |\beta_p| = 1$.

The symmetric square $L$-function $L(s,\symtwo f)$ is defined by
$$
L(s, \symtwo f) = \zeta(2s) \sum_{n=1}^{\infty} \frac{\lambda_f(n^2)}{n^s} = \prod_p \Big(1 - \frac{\alpha_p^2}{p^s}\Big)^{-1}\Big(1 - \frac{1}{p^s}\Big)^{-1}\Big(1 - \frac{\beta_p^2}{p^s}\Big)^{-1}.
$$  
From the work of Shimura \cite{Shi_2} we know that $L(s,\symtwo f)$ has an analytic continuation to the entire complex plane, and 
satisfies a functional equation connecting $s$ and $1-s$: namely, with $\Gamma_{\Bbb R}(s) = \pi^{-s/2} \Gamma(s/2)$, 
$$ 
\Lambda(s,\symtwo f) = \Gamma_{\Bbb R} (s+1) \Gamma_{\Bbb R}(s+k-1) \Gamma_{\Bbb R}(s+k) L(s,\symtwo f) = \Lambda(1-s,\symtwo f).
$$ 
  Moreover, Gelbart and Jacquet  \cite{GelJac} have shown that $L(s,\symtwo f)$ arises as the $L$-function of a cuspidal automorphic representation of $GL(3)$.    Invoking the Rankin-Selberg $L$-function attached to $\symtwo f$, a standard argument establishes the 
classical zero-free region for $L(s,\symtwo f)$, with the possible exception of a real Landau-Siegel zero (see Theorem 5.42 of \cite{IK}). 
The work of Hoffstein and  Lockhart \cite{HofLoc} (especially the appendix by Goldfeld, Hoffstein and Lieman) has ruled out the existence of Landau-Siegel zeroes for this family. 
Thus,  for a suitable constant $c > 0$, the region
$$
\CR = \Big\{s = \sigma + it : \sigma \ge 1 - \frac{c}{\log k(1+|t|)}\Big\}
$$
does not contain any zeroes of $L(s,\symtwo f)$ for any Hecke eigenform $f$ of weight $k$.

Lastly, we shall need a ``log-free'' zero-density estimate for this family, which follows from the work of Kowalski and Michel (see \cite{KM}, and also the 
recent works of Lemke Oliver and Thorner \cite{OliTho}, and Motohashi \cite{M}).

\begin{lemma} 
\label{lem2}  There exist absolute constants $B$, $C$, and $c$ such that for all $1/2 \le \alpha \le 1$, and any $T$  we have 
$$ 
| \{ \rho = \beta+i\gamma : \ \ L(\rho, \symtwo f ) = 0, \ \beta \ge \alpha, \ |\gamma |\le T \} | \le C (T+1)^B k^{c(1-\alpha)}. 
$$ 
\end{lemma}

The special value $L(1,\symtwo f)$ shows up naturally when comparing the $L^2$ normalization and 
Hecke normalization of a modular form.  Suppose $f$ has been normalized in such a way that 
$$ 
\int_{X} y^k |f(z)|^2 \frac{dx \ dy}{y^2} = 1. 
$$ 
Then the Fourier expansion of $f(z)$ is given by (see, for example, Chapter 13 of \cite{Iwaniec}) 
\begin{equation}
\label{eq:fourier_expansion}
f(z) = C_f \sum_{n=1}^{\infty} \lambda_f(n) (4 \pi n)^{\frac{k-1}{2}} e(nz),
\end{equation}
where
$$
C_f = \Big(\frac{2\pi^2}{\Gamma(k)L(1,\symtwo f)}\Big)^{1/2}.
$$
We can now state our main lemma, which refines Lemma 2 of \cite{HolSou}, and allows us to estimate $L(1,\symtwo f)$ by a 
suitable Euler product.   Below we use the notation $g \asymp h$ to denote $g \ll h$ and $h\ll g$.  
 
 \begin{lemma}
\label{lem3}
For any Hecke eigenform $f$ of weight $k$ for the full modular group, we have 
$$
L(1,\symtwo f) \asymp \exp \Big(\sum_{p \leq k} \frac{\lambda_f (p^2)}{p}\Big).
$$
 \end{lemma}
Recall that $g \asymp h$ means $g \ll h$ and $h \ll g$.
\begin{proof}
Let $1 \leq \sigma \leq \frac{5}{4}$, and consider for some $c > 0$ and $x \ge 1$, the integral
\begin{equation}
\label{eq:contour_log_L}
\frac{1}{2\pi i} \int_{c - \infty}^{c + \infty} - \frac{L'}{L} (s + \sigma, \symtwo f) (s+1)\Gamma(s) x^s ds, 
\end{equation}
which we shall evaluate in two ways.  Here we shall take $x=k^A$ for a suitably large constant $A$.   
On one hand, we write
\[
- \frac{L'}{L}(s, \symtwo f) = \sum_{n=1}^{\infty} \frac{\Lambda_{\symtwo f}(n)}{n^s}
\]
where $\Lambda_{\symtwo f}(n) = 0$ unless $n = p^k$ is a prime power, in which case 
$$ 
\Lambda_{\symtwo f}(p^k) = (\alpha_p^{2k} + 1 +\beta_p^{2k}) \log p, 
$$ 
so that $|\Lambda_{\symtwo f}(n)| \le 3 \Lambda(n)$ for all $n$.   Using this in \eqref{eq:contour_log_L}, and integrating term by term, using 
$$ 
\frac{1}{2\pi i} \int_{(c)} (s+1) \Gamma(s) y^s ds = e^{-1/y} \Big( 1+ \frac 1y \Big),
$$ 
we obtain 
\begin{align}
\label{eq:contour_log_L_right}
    \frac{1}{2\pi i}\int_{(c)} - \frac{L'}{L}(s + \sigma, \symtwo f) (s+1)\Gamma(s) x^s ds &= \sum_{n=2}^{\infty} \frac{\Lambda_{\symtwo f}(n)}{n^{\sigma}} e^{-n/x}\Big(1+\frac{n}{x}\Big). 
 \end{align}
 
On the other hand, shift the line of integration in \eqref{eq:contour_log_L} to $\Re(s) = -3/2$. We encounter poles at $s = 0$, and at $s = \rho - \sigma$ for non-trivial zeroes $\rho = \beta + i\gamma$ of $L(s, \symtwo f)$. Computing these residues, we see that \eqref{eq:contour_log_L} equals
\begin{equation} 
\label{4} 
- \frac{L'}{L}(\sigma, \symtwo f) - \sum_{\rho} x^{\rho - \sigma}(\rho-\sigma+1) \Gamma(\rho - \sigma) + \frac{1}{2\pi i } \int_{(-3/2)} -\frac{L'}{L}(s + \sigma, \symtwo f) x^s (s+1) \Gamma(s) ds.
\end{equation} 
Differentiate the functional equation of $L(s, \symtwo f)$ logarithmically, and use Stirling's formula. Thus with $s = -\frac{3}{2} + it$ we obtain that
$$
-\frac{L'}{L}(s + \sigma, \symtwo f) \ll \log (k(1+|t|) + \left|\frac{L'}{L} ( 1 - s - \sigma, \symtwo f)\right| \ll \log (k(1+|t|)).
$$
Therefore the integral in \eqref{4} may be bounded by $O((\log k) x^{-3/2})$, and we conclude that
\begin{equation}
\label{5}
\sum_{n} \frac{\Lambda_{\symtwo f}(n)}{n^{\sigma}} e^{-n/x}\Big(1+\frac nx\Big) = - \frac{L'}{L}(\sigma, \symtwo f) - \sum_{\rho} x^{\rho - \sigma} (\rho+1-\sigma) \Gamma(\rho - \sigma) + O(x^{-3/2} \log k).
\end{equation}

We now bound the sum over zeros in \eqref{5}.  Write $\rho = \beta + i\gamma$, and split into terms with $n \le |\gamma | < n+1$, where $n=0$, $1$, $2$, $\ldots$.   
If $n \le |\gamma| <n+1$, we may check using the exponential decay of the $\Gamma$-function that 
$$ 
|\rho-\sigma +1 | |\Gamma(\rho-\sigma)| \ll (\sigma -\beta)^{-1} e^{-n}.   
$$ 
Therefore the contribution of zeros from this interval is 
$$ 
\ll \sum_{n \le |\gamma| < n+1} \frac{x^{\beta- \sigma}}{\sigma-\beta} e^{-n}.
$$ 
Splitting the zeros further based on $1-(j+1)/\log k \le \beta < 1- j/\log k$ (and using the zero free region, so that $\sigma-\beta \gg (j+1)/\log k$) the 
above is 
$$ 
 \ll e^{-n} \sum_{j=0}^{\log k} \frac{x^{1-\sigma-j/\log k}}{(j+1)/\log k} 
|\{ \beta+i\gamma: \ 1-(j+1)/\log k \le \beta < 1-j/\log k, \ n \le |\gamma| < n+1\}. 
$$ 
Now using the log-free zero density estimate from Lemma \ref{lem2}, and recalling that $x=k^A$, the quantity above is 
$$ 
\ll e^{-n} x^{1-\sigma} \log k \sum_{j=0}^{\log k}\frac{e^{-jA}}{j+1} (n+1)^B k^{c (j+1)/\log k}  \ll (n+1)^B e^{-n} x^{1-\sigma} \log k,
$$ 
provided $A \ge c+ 1$ is large enough.  Now summing over $n$, we conclude that the sum over zeros in \eqref{5} is 
$\ll x^{1-\sigma} \log k$.  

Use this bound in \eqref{5}, and integrate that expression over $1\le \sigma \le 5/4$.  It follows that 
$$ 
\log L(1,\symtwo f) = \sum_{n =2}^{\infty} \frac{\Lambda_{\symtwo f}(n)}{n\log n} e^{-n/x} \Big(1+ \frac nx\Big) +O(1) = 
\sum_{p\le x} \frac{\lambda_f (p^2)}{p} + O(1),  
$$ 
since the contribution of prime powers above is easily seen to be $O(1)$, and since 
$$ 
\sum_{p\le x} \frac{1}{p} \Big| 1 - e^{-p/x} \Big(1+\frac px\Big) \Big| + \sum_{p >x} \frac 1p e^{-p/x} \Big(1+\frac px\Big) = O(1). 
$$ 
Exponentiating, we obtain  
$$ 
L(1,\symtwo f) \asymp \exp\Big( \sum_{p\le x} \frac{\lambda_f(p^2)}{p} \Big) \asymp \exp\Big( \sum_{p\le k} \frac{\lambda_f(p^2)}{p} \Big), 
$$ 
since $x=k^A$, and $\sum_{k< p\le k^A} 1/p \ll 1$.   This concludes our proof.   
\end{proof}

\section{Proof of Theorem \ref{thm:main-thm}}
\noindent The Fourier expansion \eqref{eq:fourier_expansion} and the Parseval formula give
\begin{align}
\label{6}
\int_0^1 y^k |f(z)|^2 dx & = \frac{C_f^2}{4\pi} \sum_{n=1}^{\infty} \frac{\lambda_f(n)^2}{n} (4\pi n y)^k e^{-4 \pi n y} \nonumber \\ 
&\ll \frac{1}{\Gamma(k) L(1,\symtwo f)} \sum_{n=1}^{\infty} \frac{\lambda_f(n)^2}{n} (4\pi ny)^{k }e^{-4\pi ny}. 
\end{align}

For $\xi \ge 0$, note that 
\begin{equation} 
\label{7} 
\frac{\xi^k e^{-\xi}}{\Gamma(k)} \asymp \sqrt{k} \Big(\frac{\xi}{k}\Big)^k e^{k-\xi} \ll \begin{cases} 
\sqrt{k} \exp(-(k-\xi)^2/(4k)) &\text{if  } \xi \le 2k \\ 
\sqrt{k} (e/2)^{k-\xi} &\text{if } \xi > 2k,  
\end{cases}
\end{equation} 
where the first bound follows because $\log (1+t) \le t- t^2/4$ for $|t|\le 1$ (with $t= (\xi-k)/k$), and the second bound from $\log (1+t) \le t\log 2$ for $t\ge 1$.  

The estimate \eqref{7} with $\xi = 4\pi ny$ shows that the sum in \eqref{6} is concentrated around values of $n$ 
with $|4\pi n y -k |$ about size $\sqrt{k}$.   To flesh this out, let us first show that the contribution to \eqref{6} from $n$ with $4\pi ny \ge 2k$ 
is negligible.  Using the second bound in \eqref{7}, such terms $n$ contribute (using that $L(1,\symtwo f) \gg (\log k)^{-1}$, which follows from Lemma 3 or \cite{HofLoc})
$$
\ll \frac{1}{L(1,\symtwo f)} \sum_{n\ge k/(2\pi y)} \frac{\lambda_f(n)^2}{n} \sqrt{k} (e/2)^{k-4\pi n y} 
\ll \sqrt{k} \log k \sum_{n\ge k/(2\pi y)} \frac{\lambda_f(n)^2}{n} \frac{1}{n} e^{-k/10} \ll e^{-k/20}. 
$$
This contribution to \eqref{6} is clearly negligible.  

It remains to handle the contribution from those $n$ with $4\pi ny \le 2k$.  Divide such $n$ into intervals of the form $j\sqrt{k} \le |4\pi ny -k| < (j+1)\sqrt{k}$, where 
$0\le j \ll \sqrt{k}$.  We use the first bound in (7) with $\xi = 4\pi ny$, and in the range $j\sqrt{k} \le |4\pi ny -k| <(j+1) \sqrt{k}$ this gives 
$$ 
\frac{1}{\Gamma(k)}\frac{(4\pi n y)^k}{n} e^{-4\pi n y} \ll \frac{\sqrt{k} e^{-j^2/4}}{n} \ll \frac{y}{\sqrt{k}} e^{-j^2/8},
$$ 
provided $y\ge 1/k$ say.   
Thus the contribution from the terms $j\sqrt{k} \le |4\pi ny -k| < (j+1)\sqrt{k}$ is 
\begin{align} 
\label{8}
&\ll 
\frac{ye^{-j^2/8}}{\sqrt{k} L(1,\symtwo f)} \sum_{j\sqrt{k} \le |4\pi ny -k| < (j+1)\sqrt{k} } \lambda_f(n)^2. 
\end{align} 
 At this stage, we appeal to a result of Shiu (see Theorem 1  of \cite{Shi_1}) bounding averages of non-negative multiplicative functions in short intervals.  
 
 \begin{lemma} \label{lem4}  Let $g$ be a non-negative multiplicative function with (i) $g(p^l) \le A^l$ for some constant $A$, and (ii) $g(n) \ll_{\epsilon} n^{\epsilon}$ for any $\epsilon > 0$. Then for any $\delta > 0$, if $x^{\delta } \le z\le x$, we have
 $$ 
 \sum_{x < n \le x+z} g(n) \ll_{A,\delta} \frac{z}{ \log x} \exp\Big( \sum_{p\le x} \frac{g(p)}{p}\Big). 
 $$
 \end{lemma}  
 
 Applying this lemma in \eqref{8},  in the range $y\le k^{1/2-\delta}$, we may bound that quantity by 
 $$ 
 \ll \frac{ye^{-j^2/8}}{\sqrt{k} L(1,\symtwo f)} \frac{\sqrt{k}} {y\log k} \exp\Big( \sum_{ p\le k} \frac{\lambda_f(p)^2}{p} \Big). 
 $$ 
 Since $\lambda_f(p)^2= \lambda_f(p^2) +1$, the above bound when combined with Lemma \ref{lem3} yields $\ll e^{-j^2/8}$, and 
 summing this over all $j$ gives $\ll 1$.  Thus we conclude that the quantity in \eqref{6} is bounded, completing the proof of our 
 theorem.  
 
  \bibliographystyle{abbrv}
\bibliography{HorRefs}

\end{document}